\documentclass[12pt]{article}
\usepackage{amssymb,latexsym,pstcol,pst-plot}

\newcommand{\ben}{\begin{enumerate}}
\newcommand{\een}{\end{enumerate}}
\newcommand{\ble}{\begin{lem}}
\newcommand{\ele}{\end{lem}}
\newcommand{\bth}{\begin{thm}}
\renewcommand{\eth}{\end{thm}}
\newcommand{\bpr}{\begin{prop}}
\newcommand{\epr}{\end{prop}}
\newcommand{\bco}{\begin{cor}}
\newcommand{\eco}{\end{cor}}
\newcommand{\bcon}{\begin{conj}}
\newcommand{\econ}{\end{conj}}
\newcommand{\bde}{\begin{defn}}
\newcommand{\ede}{\end{defn}}
\newcommand{\bex}{\begin{exa}}
\newcommand{\eex}{\end{exa}}
\newcommand{\barr}{\begin{array}}
\newcommand{\earr}{\end{array}}
\newcommand{\btab}{\begin{tabular}}
\newcommand{\etab}{\end{tabular}}
\newcommand{\beq}{\begin{equation}}
\newcommand{\eeq}{\end{equation}}
\newcommand{\bea}{\begin{eqnarray*}}
\newcommand{\eea}{\end{eqnarray*}}
\newcommand{\bce}{\begin{center}}
\newcommand{\ece}{\end{center}}
\newcommand{\bpi}{\begin{picture}}
\newcommand{\epi}{\end{picture}}
\newcommand{\bfi}{\begin{figure} \begin{center}}
\newcommand{\efi}{\end{center} \end{figure}}
\newcommand{\capt}{\caption}
\newcommand{\bsl}{\begin{slide}{}}
\newcommand{\esl}{\end{slide}}

\newcommand{\bib}{thebibliography}
\newcommand{\pf}{{\bf Proof}\hspace{7pt}}

\newcommand{\Qqed}{\qquad\rule{1ex}{1ex}\medskip}

\newcommand{\ol}{\overline}

\newcommand{\hs}[1]{\hspace{#1}}
\newcommand{\hso}[1]{\hspace{-1pt}}

\newcommand{\case}[4]{\left\{\barr{ll}#1&\mbox{#2}\\#3&\mbox{#4}\earr\right.}

\def\<{\langle}
\def\>{\rangle}

\newcommand{\bbN}{{\mathbb N}}

\newcommand{\Tb}{\ol{T}}

\newcommand{\pib}{\ol{\pi}}

\newcommand{\dil}{\displaystyle}

\newcommand{\dm}{Discrete Math.\/}

\newcommand{\Gaa}{\put(0,0){\circle*{3}}}
\newcommand{\Gba}{\put(10,0){\circle*{3}}}
\newcommand{\Gca}{\put(20,0){\circle*{3}}}
\newcommand{\Gda}{\put(30,0){\circle*{3}}}
\newcommand{\Gea}{\put(40,0){\circle*{3}}}
\newcommand{\Gfa}{\put(50,0){\circle*{3}}}
\newcommand{\Gga}{\put(60,0){\circle*{3}}}

\newcommand{\Gac}{\put(0,20){\circle*{3}}}
\newcommand{\Gbc}{\put(10,20){\circle*{3}}}
\newcommand{\Gcc}{\put(20,20){\circle*{3}}}
\newcommand{\Gdc}{\put(30,20){\circle*{3}}}
\newcommand{\Gec}{\put(40,20){\circle*{3}}}

\newcommand{\Ggc}{\put(60,20){\circle*{3}}}

\newcommand{\Gad}{\put(0,30){\circle*{3}}}
\newcommand{\Gbd}{\put(10,30){\circle*{3}}}

\newcommand{\Gdd}{\put(30,30){\circle*{3}}}

\newcommand{\Gfd}{\put(50,30){\circle*{3}}}
\newcommand{\Ggd}{\put(60,30){\circle*{3}}}

\newcommand{\Gbe}{\put(10,40){\circle*{3}}}
\newcommand{\Gce}{\put(20,40){\circle*{3}}}

\newcommand{\Gee}{\put(40,40){\circle*{3}}}
\newcommand{\Gfe}{\put(50,40){\circle*{3}}}

\newcommand{\Gbg}{\put(10,60){\circle*{3}}}

\newcommand{\Gdg}{\put(30,60){\circle*{3}}}

\newcommand{\Gfg}{\put(50,60){\circle*{3}}}

\newcommand{\Gaabc}{\put(0,0){\line(1,2){10}}}

\newcommand{\Gaabd}{\put(0,0){\line(1,3){10}}}

\newcommand{\Gbabd}{\put(10,0){\line(0,1){30}}}

\newcommand{\Gcabc}{\put(20,0){\line(-1,2){10}}}

\newcommand{\Gcabd}{\put(20,0){\line(-1,3){10}}}

\newcommand{\Gcadd}{\put(20,0){\line(1,3){10}}}

\newcommand{\Gdaad}{\put(30,0){\line(-1,1){30}}}

\newcommand{\Gdadd}{\put(30,0){\line(0,1){30}}}

\newcommand{\Gdagd}{\put(30,0){\line(1,1){30}}}

\newcommand{\Gdabg}{\put(30,0){\line(-1,3){20}}}

\newcommand{\Gdafg}{\put(30,0){\line(1,3){20}}}

\newcommand{\Geadd}{\put(40,0){\line(-1,3){10}}}

\newcommand{\Geafd}{\put(40,0){\line(1,3){10}}}

\newcommand{\Gfafd}{\put(50,0){\line(0,1){30}}}

\newcommand{\Ggafd}{\put(60,0){\line(-1,3){10}}}

\newcommand{\Gacbe}{\put(0,20){\line(1,2){10}}}

\newcommand{\Gbcce}{\put(10,20){\line(1,2){10}}}

\newcommand{\Gccbe}{\put(20,20){\line(-1,2){10}}}

\newcommand{\Gdcce}{\put(30,20){\line(-1,2){10}}}

\newcommand{\Gecfe}{\put(40,20){\line(1,2){10}}}

\newcommand{\Ggcfe}{\put(60,20){\line(-1,2){10}}}

\newcommand{\Gadbg}{\put(0,30){\line(1,3){10}}}

\newcommand{\Gbddg}{\put(10,30){\line(2,3){20}}}

\newcommand{\Gdddg}{\put(30,30){\line(0,1){30}}}

\newcommand{\Gfddg}{\put(50,30){\line(-2,3){20}}}

\newcommand{\Ggdbg}{\put(60,30){\line(-5,3){50}}}

\newcommand{\Ggdfg}{\put(60,30){\line(-1,3){10}}}

\newcommand{\Gbedg}{\put(10,40){\line(1,1){20}}}

\newcommand{\Gcedg}{\put(20,40){\line(1,2){10}}}

\newcommand{\Geedg}{\put(40,40){\line(-1,2){10}}}

\newcommand{\Gfedg}{\put(50,40){\line(-1,1){20}}}

\newcommand{\Gbgfg}{\put(10,60){\line(1,0){40}}}

\newtheorem{thm}{Theorem}[section]
\newtheorem{prop}[thm]{Proposition}
\newtheorem{cor}[thm]{Corollary}
\newtheorem{lem}[thm]{Lemma}
\newtheorem{conj}[thm]{Conjecture}
\newtheorem{exa}[thm]{Example}

\begin{document}
\pagestyle{plain}

\title{Proper Partitions of a Polygon and $k$-Catalan Numbers
}
\author{
Bruce E. Sagan\\[-5pt]
\small Department of Mathematics\\[-5pt] 
\small Michigan State University\\[-5pt]
\small East Lansing, MI 48824-1027\\[-5pt] 
\small USA\\[-5pt]
\small \texttt{sagan@math.msu.edu}
}

\date{\today}
\maketitle

\begin{abstract}
Let $P$ be a polygon whose vertices have been colored (labeled)
cyclically with the numbers $1,2,\ldots,c$.
Motivated by conjectures of Propp, we are led to consider partitions
of $P$ into $k$-gons which are proper in
the sense that each
$k$-gon contains all $c$ colors on its vertices.  Counting the number
of proper partitions involves a generalization of the $k$-Catalan
numbers.  We also show 
that in certain cases, any proper partition can be obtained from
another by a sequence of moves called flips.  
\end{abstract}

\section{Introduction}
\label{i}

Let $\bbN$ denote the nonnegative integers.
In September of 2003, James Propp~\cite{pro:dom} proposed a series of related
problems to the Domino List, an email group discussing matters related
to tiling.   One of the problems was as follows.
\bcon
\label{propp}
Suppose the vertices of a convex polygon $P$ are labeled cyclically 
$1,2,1,2,\ldots$  Call a triangulation of $P$ proper if no triangle
is monochromatic and let $a_N$ be the number of such triangulations if
$P$ has $N+2$ vertices.  Then
$$
a_N=
\case{\dil\frac{2^n}{2n+1}{3n\choose n}}
{if $N=2n$ where $n\in\bbN$,}
{\dil\frac{2^{n+1}}{2n+2}{3n+1\choose n}}
{if $N=2n+1$ where $n\in\bbN$.\rule{0pt}{30pt}}
$$
\econ
Note that these counts are closely connected with the 
{\it  $k$-Catalan numbers\/} defined by
$$
C_{n,k}=\frac{1}{(k-1)n+1}{kn\choose n}
$$
for $n,k\in\bbN$.  The ordinary Catalan numbers are obtained when
$k=2$.  More information about $C_{n,k}$ can be found in Stanley's
text~\cite[pp.\ 168-173]{sta:ec2}.

We will prove Propp's conjectures below.  We will also generalize them
to partitions of $P$ involving $k$-gons for $k\ge4$.  First, however,
we need some terminology.  Let $P$ be a convex polygon whose
vertices have 
been colored (labeled) counterclockwise with the sequence
$1,2,\ldots,c,1,2,\ldots,c,\ldots$  We will always draw $P$
with a horizontal edge at the top and start the coloring with the left
endpoint of that edge.

A {\it partition of $P$} is the graph $\pi$ obtained by drawing some
straight line segments (chords) between vertices of $P$ in  a plane
fashion, i.e., so that no two chords intersect in $P$'s interior.  If
all the bounded regions of this graph are $k$-gons then it will be
called a {\it $k$-partition}.  A $3$-partition will be referred
to as a {\it triangulation}.  If a $k$-gon contains the top edge, then
its {\it standard reading\/} will be the sequence of its vertices read
counterclockwise starting with the left vertex of the top edge.

A $k$-partition is {\it proper} if
each $k$-gon contains all of the $c$ colors among its vertices.  In
the case of a triangulation with two colors, this means that no triangle is
monochromatic.  Two triangulations of a pentagon are shown in
Figure~\ref{tri}.  The one on the left is proper but the one on the
right is not.

\thicklines
\setlength{\unitlength}{2pt}
\bfi
\bpi(70,70)(-10,-10)
\put(30,-5){\makebox(0,0){$1$}}
\Gda 
\put(-5,30){\makebox(0,0){$2$}}
\Gad 
\put(65,30){\makebox(0,0){$2$}}
\Ggd 
\put(10,65){\makebox(0,0){$1$}}
\Gbg
\put(50,65){\makebox(0,0){$1$}}
\Gfg
\Gdaad \Gdagd \Gadbg \Ggdfg \Gbgfg
\Gdabg \Ggdbg
\epi
\hspace{40pt}
\bpi(70,70)(-10,-10)
\put(30,-5){\makebox(0,0){$1$}}
\Gda 
\put(-5,30){\makebox(0,0){$2$}}
\Gad 
\put(65,30){\makebox(0,0){$2$}}
\Ggd 
\put(10,65){\makebox(0,0){$1$}}
\Gbg
\put(50,65){\makebox(0,0){$1$}}
\Gfg
\Gdaad \Gdagd \Gadbg \Ggdfg \Gbgfg
\Gdabg \Gdafg
\epi
\capt{Two triangulations}\label{tri}
\efi

In the next section we will prove Propp's triangulation conjectures.
In fact, in all cases we will  give two proofs.
One will involve generating functions and the Lagrange Inversion
Formula~\cite[Section 5.4]{sta:ec2}.  The other will be combinatorial,
using objects counted by a generalization of the $k$-Catalan
numbers.  In Section~\ref{pk}, we will derive analogous formulae for
partitions  of $P$ into $k$-gons for $k\ge4$.

Section~\ref{f} concerns flips.  Suppose that two
$k$-gons in a partition $\pi$ share an edge so that their union is a
$(2k-2)$-gon, $Q$.  Then a partition $\pib$ is connected to $\pi$ by
a {\it flip}, written $\pi\sim\pib$, if it agrees with $\pi$ everywhere except that the
chord of $Q$ has been replaced by another chord
connecting two opposite vertices of $Q$.  The two triangulations in
Figure~\ref{tri} are connected by a flip where $Q$ is the
quadrilateral with standard reading $1,1,2,1$.
We will say that $\pi$ and
$\pib$ are {\it connected by a sequence of flips\/} if there is a
sequence $\pi=\pi_0\sim\pi_1\sim\ldots\sim\pi_l=\pib$. A
well-known theorem of K. Wagner~\cite{wag:bv} states that any two
triangulations of a polygon (uncolored) are connected by flips.  In
fact, Wagner's theorem applies to the more general case where one
allows the set of vertices of the triangulations to include points
interior to the polygon.  In Section~\ref{f} we  show that any two
$k$-partitions of $P$ are connected by a sequence of flips.  However,
if we insist that all the partitions in the sequence be proper, called
a {\it proper flip sequence}, only
triangulations with two colors can necessarily be connected.
This answers a question of
Propp~\cite{pro:dom}.
We should note that D. Thurston~\cite{thu:fdh} has considered flips of
two hexagons sharing {\it two} edges which is equivalent to flipping a
pair of dominos in a domino tiling.

The final section is devoted to comments and open questions.

\section{Triangulations}
\label{t}

We  first prove Conjecture~\ref{propp} which we restate here for
convenience.
\bth
\label{aN}
Let $a_N$ be the number of proper triangulations of an $(N+2)$-gon,
$P,$ whose vertices have been colored cyclically with 1 and 2.
Then
$$
a_N=
\case{\dil\frac{2^n}{2n+1}{3n\choose n}}
{if $N=2n$,}
{\dil\frac{2^{n+1}}{2n+2}{3n+1\choose n}}
{if $N=2n+1$.\rule{0pt}{30pt}}
$$
\eth
\pf
We consider a single edge as a proper partition of itself so $a_0=1$.
Now suppose $N=2n+1>0$ and consider 
a proper triangulation $\pi$ of $P$. 
The top edge of $P$ is labeled $11$.  So for $\pi$ to be proper, that
edge must be in a triangle with one of the vertices labeled 2.  Say
this is the $i$th 2 in the standard reading of $P$, where
$i\ge0$ (so we start numbering with zero).  Then
the two sides of the triangle split $P$ into a $(2i+2)$-gon and a
$(2n-2i+2)$-gon which are properly triangulated by $\pi$.  This gives
us the recursion
$$
a_{2n+1}=\sum_{i=0}^n a_{2i}a_{2n-2i}.
$$
Similarly if $N=2n>0$ we get
$$
a_{2n}=\sum_{i=0}^{2n-1} a_i a_{2n-1-i}.
$$

Let $x$ be a variable and consider the generating functions
$$
\barr{l}
\dil A_0=A_0(x)=\sum_{n\ge1} a_{2n} x^n,\\[20pt]
\dil A_1=A_1(x)=\sum_{n\ge0} a_{2n+1} x^n.\\
\earr
$$
Converting the two recursions into generating function equations gives
$$
\barr{l}
A_0=2x(1+A_0)A_1,\\[5pt]
A_1=(1+A_0)^2.
\earr
$$
Plugging the second equation into the first we obtain $A_0=2x(1+A_0)^3$
which is easy to solve by Lagrange Inversion.  We use the notation
$[x^n]A(x)$ for the coefficient of $x^n$ in the generating function
$A(x)$.  Then, for $n\ge1$, we get
$$
a_{2n}=[x^n] A_0=\frac{1}{n}[x^{n-1}] 2^n(1+x)^{3n} =
\frac{2^n}{n}{3n\choose n-1}
$$
which is equivalent to the first formula in the statement of the theorem.
Similarly, we can now use Lagrange Inversion on the formula for $A_1$
in terms of $A_0$ to obtain
$$
a_{2n+1}=[x^n] A_1 =\frac{1}{n}[x^{n-1}] 2^n(1+x)^{3n}\cdot 2(1+x) =
\frac{2^{n+1}}{n}{3n+1\choose n-1}
$$
which again can be manipulated into the form given above.\Qqed

When there are three colors, one can also compute the number of proper
triangulations.  However, if the number of vertices of $P$ is
congruent to one modulo three, then the cyclical labeling will result
in the top edge being labeled $11$ and so there can be no proper
triangulations.  So in that case, we modify the labeling so that the
last vertex in the standard reading of $P$ is labeled 2.  
The proof of the next result is so similar to
the one just given, we omit it.
\bth
\label{bN}
Let $b_N$ be the number of proper triangulations of an $(N+2)$-gon,
$P$, whose vertices have been colored cyclically with 1, 2, and 3
(with the last vertex colored 2 if $N+2$ is congruent to one modulo 3).
Then
$$
b_N=
\left\{\barr{cl}
\dil\frac{1}{3n+1}{4n\choose n}&\mbox{if $N=3n$,}\\
\dil\frac{2}{3n+2}{4n+1\choose n}&\mbox{if $N=3n+1$,\rule{0pt}{30pt}}\\
\dil\frac{3}{3n+3}{4n+2\choose n}&\mbox{if $N=3n+2$.\rule{0pt}{30pt}\Qqed}
\earr\right.
$$
\eth

We would now like to give combinatorial proofs of  these
results.  To do this, we recall one of the standard combinatorial
interpretations of the $k$-Catalan numbers.  If $P$ is a polygon with
$N+2$ uncolored vertices then $C_{n,k}$ is just the number of
partitions of $P$ into $n$ polygons each having $k+1$ vertices
provided  such a partition is
possible, i.e., when $N=n(k-1)$.  We now show that certain uncolored partitions
are related to proper partitions.  
(Trivially, uncolored partitions are just proper partitions with only
one color, but we seek something more substantial.)
This proof in the case 
$k=3$ was discovered independently by 
Yuliy Baryshnikov (as communicated by Propp~\cite{pro:dom}).

\bth
\label{a2n}
We have
\bea
a_{2n}&=& 2^n C_{n,3},\\
b_{3n}&=& C_{n,4}.
\eea
\eth
\pf
Of course these results follow immediately from the previous two
theorems, but we wish to give a combinatorial proof.

First consider the statement abou $a_{2n}$.  It suffices to give a
$2^n$-to-1 map from proper triangulations $\pi$ of a 2-colored $N$-gon
$P$, where $N=2n+2$, 
to partitions of $P$ into quadrilaterals.  Since $\pi$ is proper, every
triangle has exactly one edge whose endpoints are the same color.  It
follows that if we remove these edges then the result is a partition $\pi'$
of $P$ into $n$ quadrilaterals.  

Now take an arbitrary $4$-partition $\pi'$ of $P$. 
To show that $\pi'$ occurs $2^n$ times in the image of our map, note that any
quadrilateral $Q$ appearing in $P$ must have the colors on its
vertices alternate.  This is because if some edge of $Q$ had both
endpoints of the same color, then that chord would cut off a
subpolygon of $P$ with an odd number of vertices and it would be
impossible to partition that part of $P$ into quadrilaterals.  It
follows that the inverse image of $\pi'$ consists of all $\pi$ which
can be obtained by adding back either of the two diagonals in each
quadrilateral.  Since there are $n$ quadrilaterals, the map is
$2^n$-to-1 as claimed. 

To obtain the formula for $b_{3n}$ we need a bijection between 
proper triangulations $\pi$ of a 3-colored $N$-gon $P$, where $N=3n+2$,
to partitions of $P$ into pentagons.  Given $\pi$, consider the
triangle $T$ containing the top edge which is colored $12$.  Then the
third vertex of $T$ must be colored 3.  Now there is a unique second
triangle $T'$ containing the $13$ edge and a unique third triangle
$T''$ containing the $23$ edge.  The union of these three triangles
forms a pentagon whose standard reading is $1,2,3,1,2$.
Furthermore, each of the subpolygons of $P$ outside this pentagon have
$3n'+2$ vertices for some $n'$ (depending on the subpolygon) and are
cyclically labeled in the same way as $P$ up to a permutation of the
colors.  It follows that we can iterate this construction to find a
partition $\pi'$ of $P$ into pentagons.

To construct the inverse map, suppose we are given a pentagon
partition $\pi'$.  Then in each pentagon $R$ of $\pi'$ will have
its vertices colored cyclically as $i,i+1,i+2,i+3,i+4$ for some 
$1\le i\le 3$ where we are adding modulo three.  It follows that there
will be a single  color $j$ which appears only once among the vertices
of $R$ and the other two colors will both appear twice.  So there is a
unique way of making a proper triangulation of $R$, namely by adding
the two chords containing the vertex colored $j$.  Doing this in each
pentagon, produces the inverse map.  \Qqed

We would also like to have noncolored analogues of the $a_N$'s and
$b_N$'s which do not correspond to $k$-Catalan numbers.  Let
$d\in\bbN$.  Let $P$ be a polygon rooted at an edge
which we will always take  to be the top edge.
A {\it $(k,d)$-partition of $P$\/} is a partition such that all the
regions are $k$-gons except for the one containing the root edge which
is a $d$-gon.  By convention if $d=2$ then, since the root edge is
the only edge containing itself, we just have an ordinary
$k$-partition of $P$.
Define the {\it $(k,d)$-Catalan number\/} to be
$$
C_{n,k,d}=\frac{d}{(k-1)n+d}{kn+d-1\choose n}. 
$$
Note that $C_{n,k,1}=C_{n,k}$.  The numbers $C_{n,3,d}$ have appeared
in the work of Brown on nonseparable planar maps~\cite{bro:enp};
Deutsch, Feretic and Noy on directed polyominoes~\cite{dfn:dcd}; and
of Noy on noncrossing trees~\cite{noy:ent}.  As far as we
know, combinatorial interpretations have not been given to the other
$C_{n,k,d}$.  

The following result generalizes the
$k$-partition interpretation of $C_{n,k}$.  Similar generalizations
can be given for  other interpretations of the $k$-Catalan numbers.
\bth
For $n\ge0$, $d\ge1$ and $k\ge2$, let $P$ be a rooted polygon with
$n(k-1)+d+1$ uncolored vertices.  Then
$$
C_{n,k,d}=
\mbox{number of $(k+1,d+1)$-partitions of $P$ into $n$ regions which
  are $k$-gons.}
$$
\eth
\pf
The proof is much like that of Theorem~\ref{aN} so we will just sketch
it. Considering the way the $(d+1)$-gon splits $P$ leads to a
recursion for $e_{n,k,d}$ which is defined to be the right side
of the above equation.  Letting
$$
\barr{l}
E_1=E_1(x)=\dil\sum_{n\ge1}e_{n(k-1)+2}x^n,\\[20pt]
E_d=E_d(x)=\dil\sum_{n\ge0}e_{n(k-1)+d+1}x^n,
\earr
$$
for $d\ge2$ we get functional equations
$$
\barr{l}
E_1=x(1+E_1)^k\\
E_d=(1+E_1)^d.
\earr
$$
Using Lagrange Inversion completes the proof.
\Qqed

Now we can give a more definitive version of Theorem~\ref{a2n}
\bth
\label{a2nd}
For $d=1,2$ we have
$$
a_{2n+d-1}=2^n C_{n,3,d}.
$$
For $d=1,2,3$ we have
$$
b_{3n+d-1}=C_{n,4,d}.
$$
\eth
\pf
As before, we are done if we appeal to our previous theorems but we
wish to give a combinatorial proof.  The proof is similar to that of
Theorem~\ref{a2n}.  The only difference for $a_{2n+1}$ is that there
are now an odd number of triangles.  So the triangle containing the
root edge is not paired with anything, becoming the triangle in the
rooted partition counted by $C_{n,3,2}$.  

The same idea works for
$b_{3n+1}$ and $C_{n,4,2}$.  In the case of $b_{3n+2}$, one notes that
the top edge is labeled $12$ so that the triangle containing it has
$13$ as a chord of $P$.  Pairing this triangle with the one on the
opposite side of the $13$ chord gives the necessary quadrilateral for
$C_{n,4,3}$.  Note that this quadrilateral must have vertices
$1,2,3,2$ in the standard reading and the remaining triangles can be
grouped in triples to form pentagons as in the proof of
Theorem~\ref{a2n}.  Now to construct the inverse,  
the labeling of $P$ forces the quadrilateral in the rooted partition
to have the standard reading just given in order for the rest of $P$
to be partitionable into pentagons.  Finally, each pentagon can be
dissected into triangles, again as in the proof of Theorem~\ref{a2n}. 
\Qqed

\section{Partitions with $k\ge4$}
\label{pk}

Throughout this section we will assume that $c=k\ge4$.  It will also
simplify notation to write the $k$-Catalan numbers as
$$
C_{n,k}=\frac{1}{n}{kn\choose n-1}.
$$
This is equivalent to the original definition except when $n=0$ in
which case the latter is not well defined.
\bth
\label{cN}
Let $c_N$ be the number of proper $k$-partitions of an $(N+2)$-gon,
$P$, whose vertices have been colored cyclically with $1,2,\ldots,k$
where $k\ge4$.  Then $c_0=1$ and for $N\ge1$
$$
c_N=
\left\{\barr{cl}
\dil\frac{1}{n}{(k-1)^2 n\choose n-1}
    &\mbox{if $N=(k-2)kn$,}\\
\dil\frac{k-1}{n}{(k-1)^2 n + (k-2)\choose n-1}
    &\mbox{if $N=(k-2)(kn+1)$,\rule{0pt}{30pt}}\\
0&\mbox{else.\rule{0pt}{30pt}}
\earr\right.
$$
\eth
\pf
There does not exist any $k$-partition of $P$ if  $k-2$ does not divide
$N$, so clearly $c_N=0$ in this case.  Thus we may assume that $M=N/(k-2)$ is
an integer.  Dividing $M$ by $k$ we can write $M=kn+r$ for some
$n\ge0$ and $0\le r<k$.  

We claim that $c_N=0$ if $r\neq 0,1$.  We prove this by induction.
Proceeding as in the proof of Theorem~\ref{aN} we have
$$
c_N=\sum_{N_1+\cdots N_{k-1}=N-(k-2)} c_{N_1}\cdots c_{N_{k-1}}.
$$
Suppose a term in the sum is nonzero, forcing $N_i$ to be divisible by
$k-2$ for $1\le i\le k-1$.  So we write $N_i/(k-2)=M_i=kn_i+r_i$ for
each $i$.  Also we may assume that $r_i=0$ or 1 for each $i$, either
by induction or by direct inspection in the base case $N=2(k-2)$.  If we
have both an $r_i=0$ and an $r_j=1$ then in the sequence
$r_1,\ldots,r_{k-1}$ we must have a zero followed by a one or
vice-versa.  But then in the $k$-gon containing the top edge, the
edges  corresponding to these two $c_{N_i}$ form a path of length two
whose endpoints have the same color because they are at a distance
which is a multiple of $k$ counterclockwise along $P$.  So the
partition is not proper, contradicting the fact that the term is
nonzero.  So the only other possibility is that $r_i=0$ for all $i$ or
$r_i=1$ for all $i$ which correspond to $r=1$ or $r=0$ since the $N_i$
sum to $N-(k-2)$.

The rest of the proof proceeds as in Theorem~\ref{aN}.  One defines
generating functions
$$
\barr{l}
\dil C_0=C_0(x)=\sum_{n\ge1} c_{(k-2)kn} x^n,\\[20pt]
\dil C_1=C_1(x)=\sum_{n\ge0} c_{(k-2)(kn+1)} x^n\\
\earr
$$
which satisfy functional equations
$$
\barr{l}
C_0=x C_1^{k-1},\\
C_1=(C_0+1)^{k-1}.
\earr
$$
Lagrange Inversion completes the proof.\Qqed

Again, we can give a combinatorial proof of the portion of the
previous theorem related to the $(k,d)$-Catalan numbers.
\bth
\label{ck-2kn}
We have
$$
\barr{l}
c_{(k-2)kn}=C_{n,(k-1)^2,1},\\
c_{(k-2)(kn+1)}=C_{n,(k-1)^2,k-1}.
\earr
$$
\eth
\pf
For the first equality, 
it suffices to find a bijection between proper $k$-partitions $\pi$ of a
polygon $P$ with $(k-2)kn+2$ vertices and uncolored partitions
$\pi'$ of $P$ into subpolygons with $(k-1)^2+1=(k-2)k+2$ vertices.
Given $\pi$, consider the $k$-gon, $Q$, containing the top edge.  From
the combinatorial part of the proof of the previous theorem, $c_N$ is
a sum of products of $c_{N_i}$ where the associated remainders satisfy
$r_i=1$ for all $i$.  It follows that the vertices of $Q$ read
counterclockwise are $1,k,k-1,\ldots,2$.  Now glue the $k$-gons
sharing an edge with $Q$ onto $Q$ to form a polygon $R$ with
$(k-2)k+2$ vertices.  Similar considerations show that $R$'s vertices
read counterclockwise will be the same as the usual color ordering
we use for polygons.  So we can remove $R$ from $P$ and iterate this
construction.  The collection of $R$'s obtained form the desired
partition $\pi'$.

To obtain the inverse map, consider a $[(k-2)k+2]$-partition $\pi'$ of $P$.
Then each subpolygon $R$ will be labeled in the usual coloring order up to
a permutation of the colors.  So there is a unique proper $k$-partition
of $R$, namely the one obtained by drawing a chord from the 1 of the
top edge to 
the first $k$ going counterclockwise, then another chord from that $k$ to the
next possible $k-1$ going in the same direction, and so forth
(assuming  for the sake of the description that the color permutation
is the identity).   Once all of the $R$'s have been partitioned in
this manner, one obtains a proper $k$-partition $\pi$ of $P$.  It is
easy to see that this is indeed the inverse, so we are done.

For the second inequality, note that the number of $k$-gons in $\pi$
will be one more than a multiple of $k$.  So we will be able to glue
them together as before except that one, the root polygon, will be
left over.  In other regards, we have essentially the same bijection.
\Qqed

\section{Flips}
\label{f}

We will first consider uncolored partitions.  It will be useful to use
one of the other combinatorial interpretations of $C_{n,k}$ in terms
of $k$-ary trees~\cite{sta:ec2}.  A {\it $k$-ary tree}, $T$, is a
rooted, plane tree where each vertex has either $k$ children or no
children.  The  former vertices are called {\it internal\/} and the
latter {\it leaves}.  The subtree $T_v$ of $T$ {\it generated\/} by a vertex
$v$ consists of $v$ and all its descendants.  If $v$ is an internal
vertex then we let $v',v'',\ldots,v^{(k)}$ be 
its children listed left to right and let
$T_v',T_v'',\ldots,T_v^{(k)}$ denote the trees the trees they
generate, respectively.  Vertex $v'$ is called the {\it first\/} or
{\it leftmost\/} child of $v$ while $v^{(k)}$ is the {\it last\/} or
{\it rightmost}.

\thicklines
\setlength{\unitlength}{2pt}
\bfi
\btab{ccc}
\bpi(70,40)(-10,-10)
\put(30,-5){\makebox(0,0){$1$}}
\Gda 
\put(-5,30){\makebox(0,0){$2$}}
\Gad 
\put(65,30){\makebox(0,0){$2$}}
\Ggd 
\put(10,65){\makebox(0,0){$1$}}
\Gbg
\put(50,65){\makebox(0,0){$1$}}
\Gfg
\Gdaad \Gdagd \Gadbg \Ggdfg \Gbgfg
\Gdabg \Ggdbg
\put(15,15){\lightgray \circle*{3}}
\put(45,15){\lightgray \circle*{3}}
\put(20,30){\lightgray \circle*{3}}
\put(5,45){\lightgray \circle*{3}}
\put(35,45){\lightgray \circle*{3}}
\put(55,45){\lightgray \circle*{3}}
\put(30,60){\lightgray \circle*{3}}
\put(30,65){\makebox(0,0){\lightgray $r$}}
\put(15,15){\lightgray \line(1,3){5}}
\put(45,15){\lightgray \line(-1,3){10}}
\put(20,30){\lightgray \line(-1,1){15}}
\put(20,30){\lightgray \line(1,1){15}}
\put(35,45){\lightgray \line(-1,3){5}}
\put(35,45){\lightgray \line(1,0){20}}
\epi
&
\hs{40pt}\raisebox{80pt}{$\mapsto$}\hs{20pt}
&
\bpi(70,40)(-10,-10)
\Gaa \Gca 
\Gbc \Gdc
\Gce \Gee
\Gdg
\put(30,65){\makebox(0,0){$r$}}
\Gaabc \Gcabc
\Gbcce \Gdcce
\Gcedg \Geedg
\epi
\\[20pt]
\bpi(70,70)(-10,-10)
\put(30,-5){\makebox(0,0){$1$}}
\Gda 
\put(-5,30){\makebox(0,0){$2$}}
\Gad 
\put(65,30){\makebox(0,0){$2$}}
\Ggd 
\put(10,65){\makebox(0,0){$1$}}
\Gbg
\put(50,65){\makebox(0,0){$1$}}
\Gfg
\Gdaad \Gdagd \Gadbg \Ggdfg \Gbgfg
\Gdabg \Gdafg
\put(15,15){\lightgray \circle*{3}}
\put(45,15){\lightgray \circle*{3}}
\put(20,30){\lightgray \circle*{3}}
\put(40,30){\lightgray \circle*{3}}
\put(5,45){\lightgray \circle*{3}}
\put(55,45){\lightgray \circle*{3}}
\put(30,60){\lightgray \circle*{3}}
\put(30,65){\makebox(0,0){\lightgray $r$}}
\put(15,15){\lightgray \line(1,3){5}}
\put(45,15){\lightgray \line(-1,3){5}}
\put(20,30){\lightgray \line(-1,1){15}}
\put(20,30){\lightgray \line(1,3){10}}
\put(40,30){\lightgray \line(-1,3){10}}
\put(40,30){\lightgray \line(1,1){15}}
\epi
&
\hs{40pt}\raisebox{80pt}{$\mapsto$}\hs{20pt}
&
\bpi(70,70)(-10,-10)
\put(30,65){\makebox(0,0){$r$}}
\Gac \Gcc \Gec \Ggc
\Gbe \Gfe
\Gdg
\Gacbe \Gccbe \Gecfe \Ggcfe
\Gbedg \Gfedg
\epi
\etab
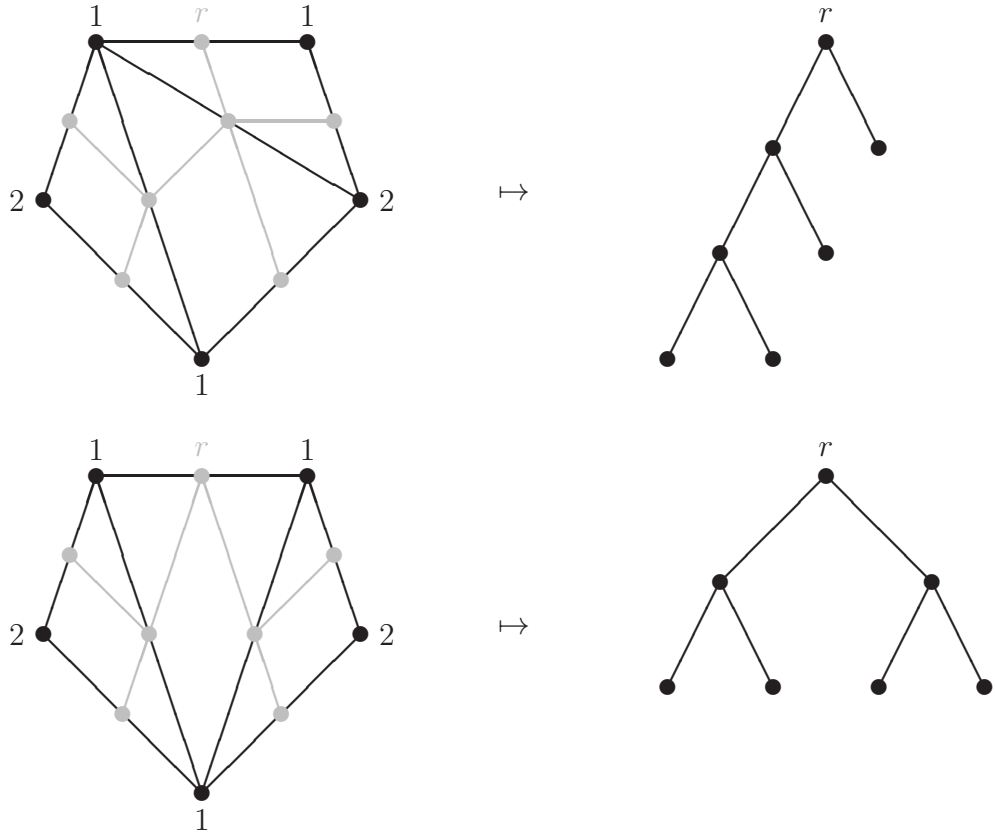
\capt{From partitions to trees}\label{trees}
\efi

It is well-known that $C_{n,k}$ counts the number of $k$-ary
trees with $n$ internal vertices.  In fact, 
there is a bijection between the partitions and trees counted by
$C_{n,k}$ which we will need.  Given at partition $\pi$ of polygon
$P$, put a tree vertex in every edge of $\pi$, including the edges of
$P$.  Now pick an edge of $P$ to contain the root vertex $r$ of $T$.  We
will always pick the top edge.   Start to build $T$ by connecting
$r$ to each of the vertices in the other edges bounding the face
containing the root edge of $P$.  This process can be iterated,
using the vertices currently adjacent to $r$ as roots of subtrees of
$T$.  An example of this construction applied to the partitions of
Figure~\ref{tri} will be found in Figure~\ref{trees}.  When the tree
is superimposed on the partition, it is shown in gray.  It is not
hard to construct the inverse for this map and thus show it is a
bijection.

We need to see what a flip does when translated into the language of
trees via this bijection.  Let $T$ be a tree and
select a vertex $v$ and one of its children $x=v^{(i)}$.  
Consider the pairwise disjoint subtrees
$$
T_v',T_v'',\ldots,T_v^{(i-1)},T_x',T_x'',\ldots,T_x^{(k)},
T_v^{(i+1)},T_v^{(i+2)},\ldots,T_v^{(k)}
$$
listed left to right in the order in which they are encountered in $T$
(i.e., in depth-first order).
Then a tree $\Tb$ is a flip of $T$, written $T\sim\Tb$ if it is
isomorphic to $T$ outside 
of $T_v$ and there is some child $y$ of $v$ such that when one makes
the list in $\Tb$ for $y$ corresponding to the above list in $T$ for
$x$, then corresponding trees in the two lists are isomorphic.  For
example, Figure~\ref{flip} shows the situation when $k=3$.  Notice
that the vertices labeled $1,2,3,4,5$ actually stand for the subtrees
generated by those vertices.

\thicklines
\setlength{\unitlength}{2pt}
\bfi
\btab{ccccc}
\bpi(50,60)(0,-10)
\put(0,-5){\makebox(0,0){$1$}}
\Gaa 
\put(10,-5){\makebox(0,0){$2$}}
\Gba 
\put(20,-5){\makebox(0,0){$3$}}
\Gca
\Gbd
\put(30,25){\makebox(0,0){$4$}}
\Gdd
\put(50,25){\makebox(0,0){$5$}}
\Gfd
\put(30,65){\makebox(0,0){$v$}}
\Gdg
\Gaabd \Gbabd \Gcabd
\Gbddg \Gdddg \Gfddg
\epi
&
\raisebox{100pt}{$\sim$}
&
\bpi(50,60)(0,-10)
\put(20,-5){\makebox(0,0){$2$}}
\Gca 
\put(30,-5){\makebox(0,0){$3$}}
\Gda 
\put(40,-5){\makebox(0,0){$4$}}
\Gea
\put(10,25){\makebox(0,0){$1$}}
\Gbd
\Gdd
\put(50,25){\makebox(0,0){$5$}}
\Gfd
\put(30,65){\makebox(0,0){$v$}}
\Gdg
\Gcadd \Gdadd \Geadd
\Gbddg \Gdddg \Gfddg
\epi
&
\raisebox{100pt}{$\sim$}
&
\bpi(50,60)(0,-10)
\put(40,-5){\makebox(0,0){$3$}}
\Gea 
\put(50,-5){\makebox(0,0){$4$}}
\Gfa 
\put(60,-5){\makebox(0,0){$5$}}
\Gga
\put(10,25){\makebox(0,0){$1$}}
\Gbd
\put(30,25){\makebox(0,0){$2$}}
\Gdd
\Gfd
\put(30,65){\makebox(0,0){$v$}}
\Gdg
\Geafd \Gfafd \Ggafd
\Gbddg \Gdddg \Gfddg
\epi
\etab
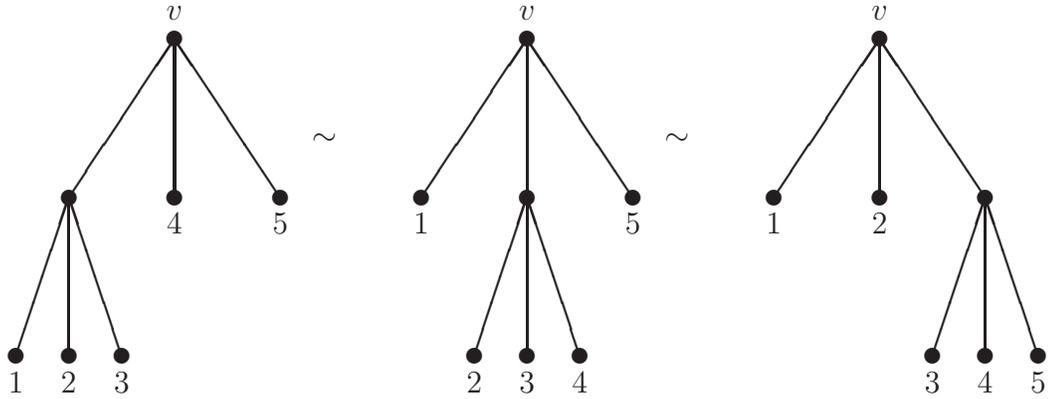
\capt{Flips when $k=3$}\label{flip}
\efi

In order to show that all $k$-ary trees with $n$ internal vertices
are connected by flips, we will need the following statistic on
trees.  The {\it left path\/}  $P$ of $T$ will be the unique path starting
at $r$ and  continuing by always taking the leftmost child.   Let
$l(T)$ denote the length of this path.  The {\it left comb}, $C$, is the
unique tree on $n$ internal vertices such that $l(C)=n$.  The first
tree in Figure~\ref{trees} is the left comb when $n=3$.  

\bth
\label{kflip}
Let $T,\Tb$ be two $k$-ary trees with $n$ internal vertices.  Then $T$
and $\Tb$ are connected by a sequence of flips.
\eth
\pf
It suffices to show that any $T$ can be connected to the left comb $C$
by a sequence of flips.  We induct on $n$.  If $n=1$ there is nothing
to prove.  Notice that
$l(T)\le l(C)$ for all $k$-ary $T$ with $n$ internal vertices, with
equality if and only if $T=C$.
So it suffices to prove that if $T\neq C$ then there is a flip such
that the resulting $\Tb$ has $l(\Tb)>l(T)$.  Since $T\neq C$ there
is some vertex $v$ on the left path of $T$ having a child $x$
such that $x\neq v'$ and $x$ is 
internal.  Using $y=v'$ for the flip creates the desired $\Tb$.
\Qqed

We will now show that when $c=2$  then any two proper
triangulations of $P$ are connect by a proper sequence of flips.  This
can be done by using the previous result and our interpretation of
colored triangulations in terms of noncolored ones.  But we prefer a
direct proof which will entail a nice characterization of the
corresponding proper trees.  Let a binary tree $T$ be {\it proper\/}
if it corresponds to a proper triangulation under the bijection
between all triangulations and all binary trees.  Then the following
result is easy to prove by induction on the number of internal nodes,
so it's proof is omitted. 
In it, $m(T)$ stands for the number of edges of $T$.
\ble
\label{properT}
A binary tree $T$ is proper if and only if for each internal vertex
$v$ either $m(T_v')$ or $m(T_v'')$ is divisible by four.
\ele

We now get a flip connection result for proper binary trees.
\bth
Let $T,\Tb$ be proper binary trees with $n$ internal nodes.  Then
there is a proper sequence of flips connecting them.
\eth
One can prove this by combining the ideas behind
Theorems~\ref{a2nd} and~\ref{kflip}.  Here we will present an
alternative direct proof.  As in the demonstration of Theorem~\ref{kflip}, it
suffices to show that given $T\neq C$ then we can connect it by a proper
sequence to some tree $U$ where $l(U)>l(T)$.  
Let $x$ and $y$ be the
right and left children of the root $r$, respectively.  By induction,
we can turn $T_x$ and $T_y$ into combs  by a proper
sequence.  Call the result $V$.  If $l(V)>l(T)$ then we are done.

Otherwise, note that $x$ is internal and $V_x''$ is a single vertex.  If
$m(V_y)$ or $m(V_x')$ is divisible by four then, by the previous lemma, 
we can apply a flip with
$v=r$ and $x,y$ playing the same roles they did in the definition
to obtain a proper tree $U$ with $l(U)>l(T)$.  If both $m(V_y)$
and $m(V_x')$ have remainder two on division by four, then do a flip
with  $x$, $x'$ and $x''$ taking the roles of $v$, $x$ and $y$,
respectively.   The resulting tree $W$ is proper and now doing the
flip with $v=r$ and $x,y$ as usual gives the desired tree $V$.
\Qqed

Connectivity by a proper sequence of flips breaks down for $c=k\ge3$.
For example, $c_{(k-2)(k+1)}$ counts the $(k-1)$-ary trees with $k+1$
internal vertices where the
root has exactly one internal child and that child has $k$ internal children.
Clearly none of these are connected by a flip.

\section{Comments and open problems}
\label{cop}

\subsection{Other labelings}

Propp~\cite{pro:dom} also conjectured a formula for the number
of proper  triangulations of a polygon
colored so that the standard reading is $m$ ones followed by $n$
twos, denoted $1^m,2^n$.  We prove it now.
\bpr
Let $d_{m,n}$ be the number of proper triangulations of a polygon $P$
colored $1^m,2^n$.  Then
$$
d_{m,n}={m+n-2\choose m-1}.
$$
\epr
\pf
If the triangle containing the top edge does not have one of the two
nodes adjacent to that edge as its third vertex, then it will split
$P$ into two parts one of which will be monochromatic making further
subdivision impossible.  This observation leads to the recursion
$d_{m,n}=d_{m-1,n}+d_{m,n-1}$ which, in conjunction with the boundary
values $d_{1,n}=d_{m,1}=1$, yields the result.
\Qqed

This raises the possibility that there may be other colorings of $P$
which will lead to nice enumerations of the corresponding proper
partitions.  One can not generalize the previous proposition directly
because for $c\ge3$ colors arranged in $c$ blocks it is easy to see that
there are no possible proper partitions.  But it would be interesting
to find other arrangements of colors which do yield nice formulae.
Note that we had to modify the cyclical labeling to get $b_{3n+2}$ to
be nonzero in Theorem~\ref{bN}.  Perhaps there are also modifications
which will do away with the zero values in Theorem~\ref{cN}.

\subsection{The case $c<k$}

The reader will have noticed that, while we permit $c<k$ in the
definition of proper, we only stated any results for this case when
$k=3$.  This is because other values lead to sequences which do not
seem to be tractable.  By way of illustration, suppose $c=3$ and
$k=4$.  Then the recursions for the corresponding sequence do not
appear to translate into simple expressions for the associated
generatiing functions.  Furthermore, the
sequence is not in Sloane's Encyclopedia of Integer
Sequences~\cite{slo:ole}.  So this avenue does not look promising.

\subsection{Other definitions}

Our definition of proper was carefully chosen to cover all cases found
so far where enumeration in closed form is possible.  But it is
conceivable that other definitions would also yield interesting
results.  For example, one might try defining proper to mean that no
$k$-gon is monochromatic.  Unfortunately, this does not seem to bear
fruit.  For example, suppose that $c=k=3$ and that $P$ is an
$(N+2)$-gon with the usual cyclic coloring.  Let
$$
b_N'=\mbox{number of triangulations of $P$ with no triangle monochromatic.}
$$
Then proceeding in the usual way using recursions,
one is led to solving the following system of generating function
equations
$$
\barr{l}
B_0=2x(1+B_0)B_2+xB_1^2,\\
B_1=(1+B_0)^2+2xB_1B_2,\\
B_2=2B_1(1+B_0).
\earr
$$
Handing the problem to Mathematica results in an output where the
solution depends on solving a quintic equation.  And the sequence
$b_N'$ is not in Sloane.

Another approach to obtaining more results would be to extend the
definition of proper to $c>k$ by saying that in this case each $k$-gon
needs to have $k$ different colors on its vertices.  We have checked
the case $c=4$ and $k=3$, but run up against the same problem as in
the previous paragraph.  However, it seems that there should be some
definition of proper which would give colored versions of all the
$(k,d)$-Catalan numbers and not just those with parameters
$((k-1)^2,1)$ or $((k-1)^2,k-1)$.

\subsection{Eliminating induction}

In the proof of Theorem~\ref{a2n} the proof that $a_{2n}=2^n C_{n,3}$
was a global construction involving flipping the diagonals of
quadrilaterals.  By contrast the proof of $b_{3n}=C_{n,4}$, while
still combinatorial, was inductive.  It would be pleasing to have a
noninductive proof of the later result.  The same applies to the
identities in Theorem~\ref{ck-2kn}. 

\subsection{Proper flip sequences}

It is disappointing that two proper trees can only be connected by a
sequence of proper flips in the case $c=2$, $k=3$.  But perhaps there
are some other simple moves which would suffice to connect 
proper trees in more cases.  The trees in the counterexample at the end of the
previous section are all connected by rotations about the root.
There are still examples where even flipping and rotation
are not enough to connect all pairs of proper trees.  But maybe a
careful analysis would lead to a small set of moves which would work.

\subsection{Tamari lattices}

One can put a partial order on the set of binary trees with a given
number of nodes by using the flips as the covering relations where $T$
is covered by $U$ if the flip taking $T$ to $U$ has $x=v'$ and $y=v''$
(in the notation of the flip definition).  These posets are in fact
lattices and have have been the object of study of a number of
authors, including Blass and Sagan~\cite{bs:mfl}, Edelman and
Reiner~\cite{er:hst}, Friedman and Tamari~\cite{ft:tfi},
Geyer~\cite{gey:tl}, Reading~\cite{rea:cl}, and Thomas~\cite{tho:tln}.
Thomas and Armstrong~\cite{tho:pc} have been looking at the analogous
structure for $k$-ary trees.

\begin{\bib}{99}

\bibitem{bs:mfl}  A. Blass and B. E. Sagan, M\"obius functions of
lattices, {\it Adv.\ in Math.} {\bf 127} (1997), 94--123.

\bibitem{bro:enp} W. G. Brown, Enumeration of non-separable planar
maps, {\it Canad.\ J. Math.\/} {\bf 15} (1963), 526--545.

\bibitem{dfn:dcd} E. Deutsch, S. Feretic, and M. Noy, Diagonally
convex directed polyominoes and even trees: a bijection and related
issues, {\it Discrete Math.\/} {\bf 180} (1998), 301--313.

\bibitem{er:hst} P. H. Edelman and V. Reiner, The higher Stasheff-Tamari
posets, {\it Mathematika} {\bf 43} (1996), 127--154.

\bibitem{ft:tfi} H. Friedman and D. Tamari, Probl\`emes
d'associativit\'e:  Une treillis finis induite par une loi
demi-associative, {\it J. Combin. Theory} {\bf 2} (1967), 215--242.

\bibitem{gey:tl} W. Geyer, On Tamari lattices, {\it \dm} {\bf 133}
(1994), 99--122.

\bibitem{noy:ent} M. Noy, Enumeration of noncrossing trees on a
circle, {\it Discrete Math.\/} {\bf 180} (1998), 301--313.

\bibitem{pro:dom} J. Propp, posting to the domino list, September, 2003.

\bibitem{rea:cl} N. Reading, Cambrian lattices, preprint,
available at {\bf http://www.arxiv.org/}, math.CO/0402086.

\bibitem{slo:ole} N. J. A. Sloane, ``The On-Line Encyclopedia of
Integer Sequences,'' available at
{\bf http://www.research.att.com/\~{\rule{1pt}{0pt}}njas/sequences/}.

\bibitem{sta:ec2} R. P. Stanley, ``Enumerative Combinatorics,
Volume 2,''  Cambridge University Press, Cambridge, 1999.

\bibitem{tho:tln} H. Thomas,  Tamari lattices and non-crossing partitions in
types $B$ and $D$, preprint,
available at {\bf http://www.arxiv.org/}, math.CO/0311334.

\bibitem{tho:pc}  H. Thomas, personal communication.

\bibitem{thu:fdh} D. P. Thurston, From dominoes to hexagons, preprint,
available at {\bf http://www.arxiv.org/}, math.CO/0405482.

\bibitem{wag:bv} K. Wagner, Bemerkungem zum Vierfarbenproblem,
{\it Jahresber.\ Deutsch.\ Math.-Verein.\/} {\bf 46} (1936), 126--132.

\end{\bib}

\end{document}